\documentclass{jglta}

\usepackage{amsmath}
\usepackage{amsthm}
\usepackage{amsfonts}
\usepackage{amssymb}
\usepackage{pstricks}
\usepackage{verbatim}
\input xy
\xyoption{all}

\JGLTAnumberwithin{equation}{section}

\allowdisplaybreaks
\newtheorem{theorem}{Theorem}[section]

\newtheorem{proposition}[theorem]{Proposition}

\theoremstyle{definition}

\allowdisplaybreaks

\newcommand{\taui}[1]{ \tau_{#1}^{i} }
\newcommand{\sectau}[1]{ \tau_{#1}^{2} }

\def\a{\alpha}
\def\b{\beta}
\def\c{\gamma}
\def\d{\delta}

\def\binl{\tau^{1}_{a}}
\def\binc{\tau^{2}_{a}}
\def\binr{\tau^{3}_{a}}

\def\alr{\tt{lr}}
\def\acr{\tt{cr}}
\def\alc{\tt{lc}}

\def\rhoA{ \rho^1 }
\def\rhoB{ \rho^2 }
\def\omegaA{ \omega^1 }
\def\omegaB{ \omega^2 }
\def\orhoA{ \overline{\rho}^{\,1} }
\def\orhoB{ \overline{\rho}^{\,2} }
\def\etaA{ \eta^1 }
\def\etaB{ \eta^2 }

\begin{document}

\renewcommand{\evenhead}{V. Abramov, R. Kerner, O. Liivapuu, and S. Shitov}
\renewcommand{\oddhead}{Algebras with ternary law of composition and their realization by cubic matrices}


\thispagestyle{empty}

\FirstPageHead{3}{*}{2008}{\pageref{firstpage}--\pageref{lastpage}}

\label{firstpage}

\Name{Algebras with ternary law of composition and their realization by cubic matrices}
\Author{V. ABRAMOV~$^a$, R. KERNER~$^b$, O. LIIVAPUU~$^a$, and S. SHITOV~$^a$}
\Address{$^a$ Institute of Mathematics, University of Tartu, Liivi 2, Tartu 50409, Estonia \\
~~E-mails: viktor.abramov@ut.ee, olgai@ut.ee, sergei.shitov@ut.ee\\
\smallskip
$^b$ LPTMC, Tour 24, Bo\^{i}te 121, 4, Place Jussieu, 75252 Paris Cedex 05,
France \\
~~E-mail: rk@ccr.jussieu.fr}

\begin{abstract}
We study partially and totally associative ternary algebras of first and second kind.
Assuming the vector space underlying a ternary algebra to be a topological space and a triple product to be continuous mapping we consider the trivial vector bundle over a ternary algebra and show that a triple product induces a structure of binary algebra in each fiber of this vector bundle. We find the sufficient and necessary condition for a ternary multiplication to induce a structure of associative binary algebra in each fiber of this vector bundle. Given two modules over the algebras with involutions we construct a ternary algebra which is used as a building block for a Lie algebra. We construct ternary algebras of cubic matrices and find four different totally associative ternary multiplications of second kind of cubic matrices. It is proved that these are the only totally associative ternary multiplications of second kind in the case of cubic matrices. We describe a ternary analog of Lie algebra of cubic matrices of second order which is based on a notion of $j$-commutator and find all commutation relations of generators of this algebra.
\end{abstract}

\section{Introduction}
A ternary algebra or triple system is a vector space $\frak A$ endowed with a ternary law of composition $\tau:{\frak A}\times{\frak A}\times{\frak A}\to {\frak A}$ which is a linear mapping with respect to each its argument, and we will call this mapping a ternary multiplication or triple product of a ternary algebra $\frak A$. Hence a ternary algebra is an algebra which closes under a suitable triple product. Obviously any binary algebra which closes under double product can be considered as a ternary algebra if one defines the ternary multiplication as twice successively applied binary one, and in this case the ternary multiplication is generated by a binary one. However there are ternary multiplications which can not be obtained as twice successively applied binary multiplication. For instance, pure imaginary numbers or elements of grading one of a superalgebra closes under triple product. A well known example of a ternary matrix algebra is the vector space $\mbox{Mat}_{m,n}$ of $m\times n$ matrices endowed with the ternary multiplication $\tau(A,B,C)=A\cdot B^{T}\cdot C$, where $A,B,C\in \mbox{Mat}_{m,n}$ and $B^T$ is transpose of the matrix $B$. Since Lie algebras play a fundamental role in physics, particular attention was given to ternary algebras when they were shown to be building blocks of ordinary Lie algebras. Given ternary algebra one can construct a Lie algebra by using the method proposed by Kantor in \cite{Kantor_1}. This method was extended to super Lie algebras in \cite{Bars-Gunadyn-1} and later was applied by the same authors in \cite{Bars-Gunadyn-2} to construct a gauge field theory by introducing fundamental fields associated with the elements of a ternary algebra.

A skew-symmetric bilinear form is an important component in the large class of algebraic structures such as Lie algebras, Grassmann algebras and Clifford algebras. For example, the Lie brackets $[\;,\;]:{\frak L}\times{\frak L}\to {\frak L}$ of a Lie algebra $\frak L$ is the skew-symmetric bilinear form, and the multiplication $(a,b)\in{\frak G}\times{\frak G}\to a\cdot b\in {\frak G}$ of a Grassmann algebra $\frak G$ restricted to the subspace of odd elements is the skew-symmetric bilinear form. A skew-symmetry of a bilinear form can be interpreted by means of the faithful representation of the symmetric group $S_2=\{e,\rho\}\to \{1,-1\}$, where $e$ is the identity permutation,  as follows: a bilinear form $\mu$ is skew-symmetric if $\mu(x_{\rho(1)},x_{\rho(2)})=(-1)\,\mu(x_1,x_2)$. Making use of this interpretation we can construct a ternary analog of a skew-symmetric bilinear form replacing $S_2$ by ${\mathbb Z}_3\subset S_3$ with its faithful representation by cubic roots of unity $j=e^{2\pi i/3}$, i.e. ${\mathbb Z}_3=\{e,\rho_1,\rho_2\}\to \{1,j,j^2\}$, where $e$ is the identity permutation and $\rho_1,\rho_2$ are the cyclic permutations, as follows: a trilinear form $\tau$ is called $j$-skew-symmetric if for any elements $a,b,c$ of a vector space $\frak A$ it satisfies
\begin{equation}
\tau(a,b,c)=j\;\tau(b,c,a)=j^2\;\tau(c,a,b)
\nonumber
\end{equation}
The notion of a $j$-skew-symmetric form can be assumed as a basis for a ternary analog of Grassmann, Clifford and Lie algebras. These ternary structures were developed in \cite{Abramov-1,Abramov-2,Kerner_3,Vainerman-Kerner} and applied to construct a ternary analog of supersymmetry algebra in
\cite{Abr-Kerner-Roy,Kerner_1,Kerner_2,LeRoy_1}.

In this paper we study algebras with ternary law of composition. In Section 2 we consider partially and totally associative ternary algebras of first and second kind. We show that a triple product of a ternary algebra induces three binary multiplications and find the sufficient and necessary condition a triple product of a ternary algebra must satisfy in order to induce the associative binary algebra. Assuming the vector space underlying a ternary algebra to be a topological space and a triple product to be continuous mapping we consider the trivial vector bundle over a ternary algebra and show that a triple product induces a structure of binary algebra in each fiber of this vector bundle. The sufficient and necessary condition a ternary multiplication must satisfy in order to induce a structure of associative binary algebra in each fiber is given in terms of the vector bundle over a ternary algebra. The relations for different kinds of partial and total associativity of a ternary algebra and induced by it binary algebras are found in the terms of the structure constants of a ternary algebra. It should be pointed out that the cohomologies of a ternary algebra of associative type are studied in \cite{Ataguema-1}.

In Section 3 we consider an algebraic structure consisting of two bimodules over unital associative algebras with involution and construct a ternary algebra by means of this algebraic structure. Choosing different modules, unital associative algebras and homomorphisms we show that this structure allows to construct a large class of ternary algebras including a ternary algebra of rectangular matrices and ternary algebras of sections of a vector bundle over a smooth finite dimensional manifold. We end the Section 3 by constructing the binary Lie algebra of matrices whose entries are the elements of bimodules and unital associative algebras. It should be mentioned that there are $n$-ary generalizations of Lie algebra which include the concepts such as $n$-ary algebra of Lie type enclosing $n$-ary Nambu algebra, $n$-ary Nambu-Lie algebra. The concept of $n$-ary Hom-algebra structure generalizing previously mentioned $n$-ary generalizations of Lie algebra is introduced and studied in \cite{Ataguema-2}. A good and detailed survey on the theory of ternary algebras can be found in \cite{Ataguema-thesis}.

It is well known that a large class of associative algebras can be constructed by means of square matrices and their multiplication. Though the rectangular matrices can be successfully used to construct a ternary algebra we think that probably more appropriate objects to construct ternary algebras are the cubic matrices. Our aim in Section 4 is to construct ternary algebras of cubic matrices and to study their structures. We find four different totally associative ternary multiplications of second kind of cubic matrices and prove that these are the only totally associative ternary multiplications of second kind in the case of cubic matrices. It is worth mentioning that our search for associative ternary multiplications of cubic matrices has shown that there is no totally associative ternary multiplication of first kind in the case of cubic matrices. I Section 5
we describe the ternary analog of Lie algebra of cubic matrices of second order by finding all commutation relations of generators of this algebra with respect to $j$-commutator.


\section{Algebras with ternary law of composition}
In this section we remind a notion of a ternary algebra and its partial or total associativity of first or second kind. Holding fixed one argument of a ternary multiplication we get the binary multiplications and study the relation between the associativity of a ternary multiplication and the associativity of induced binary multiplication. We propose to use a vector bundle approach to describe the family of binary algebras induced by a ternary algebra.

Let ${\frak A}, {\frak B}$ be complex vector spaces, and $\tau: {\frak A}\times {\frak A}\times {\frak A}\to {\frak B}$ be a $\frak B$-valued trilinear form.
We will call $\tau$ a ternary law of composition or ternary multiplication on $\frak A$ if $\tau$ is a $\frak A$-valued trilinear form.
The pair $({\frak A}, \tau)$ is said to be a ternary algebra or triple system if $\frak A$ is a complex vector space, and $\tau: {\frak A}\times{\frak A}\times{\frak A}\to {\frak A}$ is a ternary law of composition on $\frak A$.  It is obvious that relation analogous to binary associativity in the case of ternary law of composition should contain at least five elements of $\frak A$. There are three different ways to apply twice a ternary multiplication $\tau$ to ordered sequence of five elements $a,b,c,d,f\in {\frak A}$ which lead us to the following relations defining a notion of partial associativity for ternary multiplication:
\begin{eqnarray}
\tau(\tau(a,b,c),d,f)=\tau(a,b,\tau(c,d,f)) \label{lr-kind}\\
\tau(\tau(a,b,c),d,f)=\tau(a,\tau(b,c,d),f) \label{lc-kind}\\
\tau(a,\tau(b,c,d),f)=\tau(a,b,\tau(c,d,f)) \label{cr-kind}
\end{eqnarray}

Hence we have three different kinds of partially associative ternary algebra $({\frak A},\tau)$ which will be called {\alr}-partially associative ternary algebra (\ref{lr-kind}), {\alc}-partially associative ternary algebra of first kind (\ref{lc-kind}) and {\acr}-partially associative ternary algebra of first kind (\ref{cr-kind}). A ternary algebra $({\frak A},\tau)$ is said to be totally associative ternary algebra of first kind if its ternary multiplication $\tau$ satisfies any two of the relations (\ref{lr-kind})--(\ref{cr-kind}). It is obvious that in the case of totally associative ternary algebra of first kind a ternary multiplication $\tau$ satisfies the relations
\begin{equation}
\tau(\tau(a,b,c),d,f)=\tau(a,\tau(b,c,d),f)=\tau(a,b,\tau(c,d,f)) \label{first-kind}
\end{equation}
where $a,b,c,d,f\in {\frak A}$. The notion of totally associative ternary algebra of first kind can be viewed as a direct ternary generalization of classical associativity $\mu(\mu(x,y),z)=\mu(x,\mu(y,z))$, where $x,y,z$ are the elements of an algebra $({\cal A},\mu)$ with a binary law of composition $\mu:{\cal A}\times {\cal A}\to  {\cal A}$, when one applies twice algebra multiplication (binary or ternary) to ordered sequence of elements of algebra successively shifting the first (interior) multiplication from left to right and setting equal obtained products.
In this sense the notion of {\alr}-partial associativity can be considered as most similar to classical associativity whereas the notion of {\alc}-partial or {\acr}-partial associativity can be defined for the first time only in the case of ternary multiplication because in the case of binary multiplication there is no central group of two elements in the middle of a sequence $x,y,z\in {\cal A}$.

Since the notion of {\alc}-partial or {\acr}-partial associativity appears for the first time in the case of ternary multiplication there is no reason to keep the requirement of fixed order of a sequence $a,b,c,d,f\in {\frak A}$ looking for a possible analogs of associativity in the case of ternary algebras. It turns out that we get a useful notion of ternary associativity giving up the requirement of fixed order of elements in a sequence $a,b,c,d,f\in {\frak A}$. This means that unlike the case of ternary associativity of first kind we not only successively shift the first (interior) multiplication inside a sequence of elements
$a,b,c,d,f\in {\frak A}$ from left to right but at the same time permute the elements $b,c,d$ in the middle of sequence. Obviously we should use non-cyclic permutation in order to get the initial order of a sequence $a,b,c,d,f$ on the second step. This reasoning leads us to the following relations:
\begin{eqnarray}
\tau(\tau(a,b,c),d,f)=\tau(a,\tau(d,c,b),f) \label{lc-secondkind}\\
\tau(a,\tau(d,c,b),f)=\tau(a,b,\tau(c,d,f))\label{cr-secondkind}
\end{eqnarray}
A ternary algebra $({\frak A}, \tau)$ is said to be {\alc}-partially associative ternary algebra of second kind if ternary multiplication $\tau$ satisfies (\ref{lc-secondkind}) and {\acr}-partially associative ternary algebra of second kind if $\tau$ satisfies (\ref{cr-secondkind}). A ternary algebra $({\frak A}, \tau)$ is said to be totally associative ternary algebra of second kind if it is {\alr}-partially associative and either {\alc}-partially associative of second kind or {\acr}-partially associative of second kind.
Hence in the case of totally associative ternary algebra of second kind we have
\begin{equation}
\tau(\tau(a,b,c),d,f)=\tau(a,\tau(d,c,b),f)=\tau(a,b,\tau(c,d,f)) \label{second-kind}
\end{equation}

 A ternary multiplication $\tau$ of ternary algebra $({\frak A},\tau)$ has three arguments $\tau(a,b,c)$, where $a,b,c\in {\frak A}$, and if we fix one of them then $\tau$ induces the binary multiplication on ${\frak A}$. It is obvious that this allows us to split the ternary multiplication $\tau$ into three binary ones. We can study the structure of ternary multiplication $\tau$ from this point of view by making use of known concepts and methods of the theory of binary algebras. Given an element $a\in {\frak A}$ a ternary multiplication $\tau$ induces three binary multiplications $\binl,\binc,\binr$ on $\frak A$ defined as follows:
\begin{equation}
\binl (b,c)=\tau(a,b,c),\quad
     \binc (b,c)=\tau(b,a,c),\quad
        \binr (b,c)=\tau(b,c,a) \label{binary_operations}
\end{equation}
where $b,c\in {\frak A}$. The binary multiplications $\tau^1_a,\tau^2_b,\tau^3_c$ are not independent because of the relations
\begin{equation}
\tau^1_a(b,c)=\tau^2_b(a,c)=\tau^3_c(a,b)
\end{equation}
A vector space $\frak A$ equipped with the binary multiplication $\taui{a}, i=1,2,3$ becomes the binary algebra which will be denoted by $({\frak A},\taui{a})$. Considering an element $a$ in $({\frak A},\taui{a})$ as a parameter ranging within a vector space $\frak A$ we have three families of binary algebras $({\frak A},\binl),({\frak A},\binc),({\frak A},\binr)$ induced by a ternary multiplication $\tau$. The family of binary algebras $({\frak A},\taui{a})$ is said to be an associative family of binary algebras induced by a ternary algebra $({\frak A},\tau)$ if for any $a,b,c,d,f\in {\frak A}$ it holds
\begin{equation}
\taui{a}(\taui{b}(c,d),f)=\taui{b}(c,\taui{a}(d,f)) \label{associativefamily}
\end{equation}
Taking $a=b$ in the previous relation we see that each associative family of binary algebras $({\frak A},\taui{a})$ is the family of associative binary algebras.

It is useful to describe the above mentioned families of binary algebras in terms of vector bundle. For this purpose we will assume that ${\frak A}$ is a topological vector space, and a ternary multiplication $\tau:{\frak A}\times{\frak A}\times{\frak A}\to {\frak A}$ is a continuous mapping. Let us consider the direct product
${\frak E}={\frak A}\times {\frak A}$ as the trivial vector bundle over the base space ${\frak A}$ with the fiber $\frak A$ and the projection $\pi:{\frak E}\to {\frak A}$ defined as usual $\pi(p)=a$, where $p=(a,b)\in {\frak E}$. Any fiber $\pi^{-1}(a)$ of ${\frak E}$ is isomorphic to ${\frak A}$, and we will denote this isomorphism at a point $a$ of the base space $\frak A$ by $\phi_a$, i.e. $\phi_a:{\pi}^{-1}(a)\to {\frak A}$ and $\phi_a(p)=b$, where $p=(a,b)\in \pi^{-1}(a)$. Let $a,b\in {\frak A}$ be two points of the base space of a vector bundle ${\frak E}$. Then \begin{equation}
\phi^a_{b}=\phi^{-1}_{b}\circ\phi_{a}:\pi^{-1}(a)\to \pi^{-1}(b)
\end{equation}
is the isomorphism between two fibers.

In order to apply the constructed vector bundle $\frak E$ to describe the families of binary algebras induced by a ternary algebra $({\frak A},\tau)$ within the framework of a single structure we assume that the base space $\frak A$ of this bundle is equipped with a ternary multiplication $\tau$.  For any point $a\in {\frak A}$ of the base space a fiber $\pi^{-1}(a)\subset {\frak E}$ at this point is endowed with one of the binary multiplications $\binl,\binc,\binr$ which we carry over from the family of binary algebras to fibers of $\frak E$ by requiring $\phi_a$ to be an isomorphism of algebras, i.e.
\begin{equation}
\phi_a(\taui{a}(p,q))=\taui{a}(\phi_a(p),\phi_a(q))
\end{equation}
where $p,q\in \pi^{-1}(a)$. If each fiber of $\frak E$ is endowed with a binary multiplication $\tau^i$ then in order to emphasize this algebraic structure of fibers we will denote the corresponding vector bundle by ${\frak E}_i$. Thus ${\frak E}_i=({\frak A},\tau)\times ({\frak A},\tau^i)$, where the base space (the first factor in the direct product) is a ternary algebra $({\frak A},\tau)$, and a fiber $\pi^{-1}(a)$ is the binary algebra $({\frak A},\taui{a})$. We will call ${\frak E}_i$, where $i=1,2,3$, the vector bundle of binary algebras over a ternary algebra $({\frak A},\tau)$. A section $\xi$ of the vector bundle ${\frak E}_i$ is a continuous mapping $\xi:{\frak A}\to {\frak E}_i$ satisfying $\pi\circ \xi=\mbox{id}_{\frak A}$, and the vector space of continuous sections will be denoted by $\Gamma({\frak E}_i)$. Evidently this vector space equipped with the binary multiplication
\begin{equation}
\taui\,(\xi,\eta)\,(a)=\taui{a}\big(\xi(a),\eta(a)\big)
\end{equation}
where $\xi,\eta\in \Gamma({\frak E}_i)$, is the binary algebra.

The notion of an associative family of binary algebras defined by (\ref{associativefamily}) can be described in the terms of vector bundle ${\frak E}_i$. Let $p,q\in {\frak E}_i$ be two points of a vector bundle ${\frak E}_i$ and $a\in {\frak A}$. A vector bundle of binary algebras ${\frak E}_i$ is said to be an associative vector bundle of binary algebras over a ternary algebra $({\frak A},\tau)$ if for any $p,q\in {\frak E}_i$ and $a\in {\frak A}$ it holds
\begin{equation}
\phi_{\pi(q)}\circ\taui{\pi(q)}\Big(\phi^{\pi(p)}_{\pi(q)}\big(\taui{\pi(p)}(p,\phi^{-1}_{\pi(p)}(a))\big),q\Big)=    \phi_{\pi(p)}\circ\taui{\pi(p)}\Big(p,\phi^{\pi(q)}_{\pi(p)}\big(\taui{\pi(q)}(\phi^{-1}_{\pi(q)}(a),q)\big)\Big) \label{associative_family}
\end{equation}
Particularly (\ref{associative_family}) implies the associativity of a fiber $\pi^{-1}(a)$ for any $a\in {\frak A}$ if we take $\pi(p)=\pi(q)$ in (\ref{associative_family}), i.e. any associative vector bundle of binary algebras ${\frak E}_i$ is a vector bundle of associative binary algebras whereas the converse is generally not true.
Now it is natural to pose a question concerning the associativity of induced binary algebras $({\frak A},\taui{a})$ provided a ternary algebra $({\frak A},\tau)$ is partially or totally associative of first or second kind.
\begin{proposition}
\label{lr-associativity}
A ternary algebra $({\frak A},\tau)$ is {\alr}-partially associative ternary algebra if and only if  ${\frak E}_2$ is the associative vector bundle of binary algebras over a ternary algebra $({\frak A},\tau)$. Particularly if a base space
$({\frak A},\tau)$ is {\alr}-partially associative ternary algebra then each fiber of the vector bundle ${\frak E}_2$ is an associative binary algebra and the binary algebra of sections of this bundle $\Gamma({\frak E}_2)$ is associative algebra.
\end{proposition}
\noindent
Indeed the left side of (\ref{associative_family}) can be transformed as follows:
\begin{eqnarray}
\phi_{\pi(q)}\circ\sectau{\pi(q)}\Big(\phi^{\pi(p)}_{\pi(q)}\big(\sectau{\pi(p)}(p,\phi^{-1}_{\pi(p)}(a))\big),q\Big)=
\phi_{\pi(q)}\circ\sectau{\pi(q)}\Big(\sectau{\pi(p)}\big(\phi^{-1}_{\pi(q)}\circ\phi_{\pi(p)}(p),\phi^{-1}_{\pi(q)}(a)     \big),q\Big)\nonumber\\
=\sectau{\pi(q)}\Big(\sectau{\pi(p)}\big(\phi_{\pi(p)}(p),a\big),\phi_{\pi(q)}(q)\Big)=
  \tau\Big(\tau\big(\phi_{\pi(p)}(p),\pi(p),a\big),\pi(q),\phi_{\pi(q)}(q)\Big) \nonumber
\end{eqnarray}
Analogously for the right side of (\ref{associative_family}) we have
\begin{eqnarray}
\phi_{\pi(p)}\circ\taui{\pi(p)}\Big(p,\phi^{\pi(q)}_{\pi(p)}\big(\taui{\pi(q)}(\phi^{-1}_{\pi(q)}(a),q)\big)\Big)=
  \tau\Big(\phi_{\pi(p)}(p),\pi(p),\tau\big(a,\pi(q),\phi_{\pi(q)}(q)\big)\Big)
\nonumber
\end{eqnarray}
and this proves the {\alr}-partial associativity of a ternary algebra $({\frak A},\tau)$.

It is well known that any associative binary algebra is the Lie algebra with respect to the commutator defined with the help of a binary multiplication of this algebra, and the associativity of a binary multiplication implies the Jacoby identity for the commutator. It follows from the Proposition \ref{lr-associativity} that if a ternary algebra $({\frak A,\tau})$ is {\alr}-partially associative ternary algebra then each fiber $\pi^{-1}(a)$, where $a\in{\frak A}$, of the vector bundle ${\frak E}_2$ is the Lie algebra with respect to the commutator $[\;,\;]_a$ defined  by
\begin{equation}
[p,q]_a=\binc(p,q)-\binc(q,p)
\end{equation}
where $p,q\in \pi^{-1}(a)$. Clearly the associative binary algebra of sections $\Gamma({\frak E}_2)$ is the Lie algebra under the commutator
\begin{equation}
[\xi,\eta](a)=[\xi(a),\eta(a)]_a$, \quad where\quad $\xi,\eta\in \Gamma({\frak E}_2)
\end{equation}

A ternary algebra $({\frak A},\tau)$ is said to be a ternary algebra of Lie type of first kind (of second kind) if for any $a_1,a_2,a_3\in {\frak A}$ it holds
\begin{equation}
\sum_{\rho\in S_3} \tau(a_{\rho(1)},a_{\rho(2)},a_{\rho(3)})=0\quad
   \Big(\sum_{\rho\in {\mathbb Z}_3} \tau(a_{\rho(1)},a_{\rho(2)},a_{\rho(3)})=0\Big)
   \label{grassmann_type}
\end{equation}
where $S_3$ is the symmetry group of third order and ${\mathbb Z}_3$ is its cyclic subgroup. Clearly any ternary algebra of Lie type of second kind is a ternary algebra of Lie type of first kind whereas the converse is generally not true. From (\ref{grassmann_type}) it follows that any element $a$ of a ternary algebra of Lie type (of first or second order) satisfies $a^3=\tau(a,a,a)=0$. It is pointed out in the Introduction that we can construct a ternary analog of the notion of skew-symmetry by means of a faithful representation of ${\mathbb Z}_3$ by cubic roots of unity. Let $j=e^{\frac{2\pi i}{3}}\in {\mathbb C}$ be the primitive cubic root of unity. A ternary multiplication $\tau$ of a ternary algebra $({\frak A},\tau)$ is said to be $j$-skew-symmetric if
\begin{equation}
\tau(a,b,c)=j\;\tau(b,c,a)=j^2\;\tau(c,a,b)
\label{cyclic j-symmetric}
\end{equation}
where $a,b,c\in {\frak A}$. If a ternary multiplication $\tau$ of $({\frak A},\tau)$ is $j$-skew-symmetric then $({\frak A},\tau)$ is a ternary algebra of Lie type of second order. Indeed in this case we have
\begin{equation}
\tau(a,b,c)+\tau(b,c,a)+\tau(c,a,b)=\tau(a,b,c)+j^2\,\tau(a,b,c)+j\,\tau(a,b,c)=0
\end{equation}
We see that the notion of $j$-skew-symmetric ternary multiplication is based on the faithful representation of the cyclic group ${\mathbb Z}_3$ by cubic roots of unity. Given a ternary algebra $({\frak A},\tau)$ we can make it the ternary algebra of Lie type of second order by endowing it with the ternary $j$-brackets or the ternary $j$-commutator which is defined by
\begin{equation}
[a,b,c]=\tau(a,b,c)+j\,\tau(b,c,a)+j^2\,\tau(c,a,b)
\end{equation}
If $\frak A$ has an involution $\ast:{\frak A}\to{\frak A}$ then $\tau$ will be called Hermitian if it satisfies $\tau(a,b,c)=\tau^{*}(c,b,a)$. Hence  a Hermitian $j$-skew-symmetric ternary multiplication $\tau$ satisfies (\ref{cyclic j-symmetric}) and
\begin{equation}
\tau(a,b,c)=\tau^{*}(c,b,a)=j^2\;\tau^{*}(b,a,c)=j\;\tau^{*}(a,c,b)
\end{equation}

Let us suppose that $({\frak A},\sigma)$ is a {\alr}-partially nonassociative ternary algebra, i.e. in general we have $\sigma(\sigma(a,b,c),d,f)\neq\sigma(a,b,\sigma(c,d,f))$, where $a,b,c,d,f\in{\frak A}$.
A ternary algebra $({\frak A},\sigma)$ is said to be a ternary algebra of Jordan type if its ternary multiplication $\sigma$ satisfies the following identities:
\begin{eqnarray}
\lefteqn{\sigma(a,b,c)=\sigma(c,b,a)} \label{first_identity} \\
\lefteqn{\sigma(\sigma(a,b,c),b,\sigma(a,b,a))=\sigma(a,b,\sigma(c,b,\sigma(a,b,a)))}\label{second_identity}
\end{eqnarray}
where $a,b,c\in {\frak A}$.
It is easy to see that if $({\frak A},\sigma)$ is a ternary algebra of Jordan type then for any $a\in {\frak A}$ the binary algebra $({\frak A},\sigma^2_{a})$ is the Jordan algebra. Indeed in this case the identities  (\ref{first_identity},\ref{second_identity}) take on the form
\begin{equation}
\sigma^2_{a}(b,c)=\sigma^2_{a}(c,b),\quad \sigma^2_{a}(\sigma^2_{a}(b,c),\sigma^2_{a}(b,b))=\sigma^2_{a}(b,\sigma^2_{a}(c,\sigma^2_{a}(b,b)))
\end{equation}
\begin{proposition}
If $({\frak A},\tau)$ is {\alr}-partially associative ternary algebra then the ternary algebra $({\frak A},\sigma)$, where $\sigma(a,b,c)=\tau(a,b,c)+\tau(c,b,a)$ is the ternary algebra of Jordan type.
\end{proposition}

We see that having fixed one variable in a triple product $\tau(a,b,c)$ of a ternary algebra $({\frak A},\tau)$ we can study the structure of a ternary multiplication $\tau$ by splitting it into three binary ones. What kind of structures induces a ternary multiplication of $({\frak A},\tau)$ if one fixes two variables in $\tau(a,b,c)$? Obviously fixing two variables we get the linear operator acting on ${\frak A}$, and this is the second way for studying the structure of a ternary multiplication. Let $\mbox{Lin}({\frak A})$ be the algebra of linear operators of the vector space $\frak A$. Given a pair $(a,b)\in {\frak A}\times {\frak A}$ we define the linear operators $L^{i}(a,b):{\frak A}\to {\frak A}$, where $i=1,2,3$, as follows:
\begin{equation}
L^{1}(a,b)\cdot c=\tau(c,a,b),\quad
     L^{2}(a,b)\cdot c=\tau(a,c,b),\quad
         L^{3}(a,b)\cdot c=\tau(a,b,c)
\end{equation}
Actually these operators are not independent because for any $a,b,c\in{\frak A}$ we have the relations
\begin{equation}
L^1(c,b)\cdot a=L^2(a,b)\cdot c=L^3(a,c)\cdot b
\end{equation}
It is easy to see that for every $i$ the linear operator $L^i(a,b)$ is bilinear with respect to its variables $a,b$, and therefore the family of linear operators $\{L^i(a,b)\}_{(a,b)\in{\frak A}\times{\frak A}}$ determines the bilinear mapping $L^{i}:{\frak A}\times {\frak A}\to \mbox{Lin}({\frak A})$. On the other hand if there is a vector space $\frak A$ equipped with a bilinear mapping $L:{\frak A}\times{\frak A}\to \mbox{Lin}({\frak A})$ then one can construct the ternary algebra $({\frak A},\tau)$ by letting $\tau(a,b,c)=L(a,b)\cdot c$.

Now we can introduce an analog of identity element for a ternary algebra $({\frak A},\tau)$ by means of the bilinear mappings $L^i:{\frak A}\times{\frak A}\to {\frak A}$, $i=1,2,3$. Indeed given an element $a\in {\frak A}$ we define ${\frak I}^{(i)}_a\subset {\frak A}\times {\frak A}$ by
\begin{equation}
{\frak I}^{(i)}_a=\{(e,\tilde e)\in {\frak A}\times {\frak A}: L^{(i)}(e,\tilde e)=\mbox{id}_{\frak A}\}
\end{equation}
where $\mbox{id}_{\frak A}\in \mbox{Lin}({\frak A})$ is the identity operator. A pair $(e,\tilde e)\in {\frak A}\times {\frak A}$ is said to be an identity $i$-pair for an element $a$ if $(e,\tilde e)\in {\frak I}^{(i)}_a$. Let ${\frak I}^{(i)}=\bigcap_{a\in {\frak A}}{\frak I}^{(i)}_a$ and ${\frak I}=\bigcap_{i=1}^3 {\frak I}^{(i)}$. If ${\frak I}^{(i)}\neq \emptyset$ then we will call an element $(e,\tilde e)\in {\frak I}^{(i)}$ an identity $i$-pair, and similarly if ${\frak I}\neq \emptyset$ then we will call an element $(e,\tilde e)\in {\frak I}$ an identity pair of a ternary algebra $({\frak A},\tau)$.

Let us now assume that the vector space $\frak A$ of a ternary algebra $(\frak{A},\tau)$
is a finite dimensional vector space, i.e.
$\frak A$ is an $r$-dimensional vector space and
$\frak{e}=\{e_1,e_2,\ldots , e_r\}$ is a basis for $\frak A.$ Then for any element $a\in\frak{ A}$ we have
$a=a^{\a}e_{\a}$, and the triple product of elements $a,b,c\in{\frak A}$ can be expressed as follows:
\begin{equation}
\tau(a,b,c)=\tau(a^{\a}e_{\a},b^{\b}e_{\b},c^{\c}e_{\c})=C^{\d}_{\a\b\c}a^{\a}b^{\b}c^{\c}e_{\d} \nonumber
\end{equation}
where $C^{\d}_{\a\b\c}$ are the structure constants of a ternary algebra $({\frak A},\tau)$ defined by
\begin{equation}
\tau(e_{\a},e_{\b},e_{\c})=C^{\d}_{\a\b\c}e_{\d}
\end{equation}
If $\frak{e}'=\{e_1',e_2',\ldots e_r'\}$ is another basis for a
vector space $\frak{A}$ and $e'_{\a}=A_{\a}^{\b}e_{\b},$ where
$A=(A_{\a}^{\b})$ is the transition matrix then
\begin{eqnarray}
\widetilde{C}_{\a\b\c}^{\d}=
A_{\a}^{\varepsilon}A_{\b}^{\zeta}A_{\c}^{\eta}\bar{A}^{\d}_{\lambda}C_{\varepsilon\zeta\eta}^{\lambda} \nonumber
\end{eqnarray}
where $\widetilde{C}_{\a\b\c}^{\d}$ are the structure constants of a ternary algebra $(\frak{A},\tau)$   with respect
to a basis $\frak{e}'$, i.e. $\tau(e_{\a}',e_{\b}',e_{c}')=\widetilde{C}_{\a\b\c}^{\d}e_{\d}'=
\widetilde{C}_{\a\b\c}^{\d}A_{\delta}^{\lambda}e_{\lambda}$, and $A^{-1}=(\bar{A}^{\d}_{\lambda})$ is the inverse matrix of $A$.

If we require a ternary algebra $({\frak A},\tau)$ to be a partially or totally associative ternary algebra either of first or second kind then this requirement leads to the relations the structure constants of $({\frak A},\tau)$ have to satisfy. These relations for different kinds of associativity of first kind have the following form:
\begin{eqnarray}
C^{\d}_{\a\b\c}C^{\nu}_{\d\lambda\mu}=C^{\nu}_{\a\b\d}C^{\d}_{\c\lambda\mu} \quad\mbox{({\alr}-partial associativity of first kind)} \nonumber\\
C^{\d}_{\a\b\c}C^{\nu}_{\d\lambda\mu}=C^{\nu}_{\a\d\mu}C^{\d}_{\b\c\lambda} \quad\mbox{({\alc}-partial associativity of first kind)} \nonumber\\
C^{\nu}_{\a\d\mu}C^{\d}_{\lambda\c\beta}=C^{\nu}_{\a\b\d}C^{\d}_{\c\lambda\mu} \quad\mbox{({\acr}-partial associativity of first kind)} \nonumber
\end{eqnarray}
It follows from the above relations that if $(\frak{A},\tau)$ is a totally associative
ternary algebra of first kind then the structure constants satisfy
\begin{gather}
C^{\d}_{\a\b\c}C^{\nu}_{\d\lambda\mu}=C^{\nu}_{\a\d\mu}C^{\d}_{\b\c\lambda}=C^{\nu}_{\a\b\d}C^{\d}_{\c\lambda\mu}
\nonumber
\end{gather}
In the case of ternary associativity of second kind we have the following relations:
\begin{eqnarray}
C^{\d}_{\a\b\c}C^{\nu}_{\d\lambda\mu}=C^{\nu}_{\a\d\mu}C^{\d}_{\lambda\c\beta} \quad\mbox{({\alc}-partial associativity of second kind)} \nonumber\\
C^{\nu}_{\a\d\mu}C^{\d}_{\lambda\c\beta}=C^{\nu}_{\a\b\d}C^{\d}_{\c\lambda\mu} \quad\mbox{({\acr}-partial associativity of second kind)} \nonumber
\end{eqnarray}
The structure constants of a totally associative ternary algebra of second kind satisfy
\begin{gather}
C^{\d}_{\a\b\c}C^{\nu}_{\d\lambda\mu}=C^{\nu}_{\a\d\mu}C^{\d}_{\lambda\c\beta}=C^{\nu}_{\a\b\d}C^{\d}_{\c\lambda\mu}
\nonumber
\end{gather}
For any $a\in {\frak A}$ a ternary algebra $({\frak A},\tau)$ induces three binary algebras $({\frak A},\taui{a}),\;i=1,2,3$, with the binary multiplications defined by relations $(\ref{binary_operations})$.
The structure constants $K_{\a\b}^{i,\c}(a)$ of binary algebra $(\frak{A},\taui{a})$ defined by
$\taui{a}(e_{\a},e_{\b})=K_{\a\b}^{i,\c}(a)e_{\c}$ can be expressed in terms of the structure constants of a ternary algebra $({\frak A},\tau)$ as follows:
\begin{equation}
K_{\a\b}^{1,\d}(a)=C_{\gamma\alpha\b}^{\d}\, a^{\gamma},\quad
   K_{\a\b}^{2,\d}(a)=C_{\alpha\gamma\b}^{\d}\, a^{\gamma},\quad
     K_{\a\b}^{3,\d}(a)=C_{\alpha\b\gamma}^{\d}\, a^{\gamma}
\end{equation}
where $a=a^{\gamma}\,e_{\gamma}$.


\section{Lie algebras from ternary algebra}

In this section we propose few different methods for constructing ternary algebras and apply these methods to construct a ternary algebra of vector fields on a smooth finite dimensional manifold and a ternary algebra of rectangular matrices. Particularly the curvature of an affine connection determines the structure of a ternary algebra on the module of vector fields on a smooth manifold. Given two modules over the algebras with involutions we construct a ternary algebra which is used to construct a Lie algebra. Our approach generalizes the approach proposed in \cite{Bars-Gunadyn-1,Bars-Gunadyn-2}, where the authors use the rectangular complex matrices.

Let $\frak A$ be a vector space over the complex numbers $\mathbb C$ and ${\frak A}^\ast$ be the dual space. Given a $\mathbb C$-multilinear mapping $T:{\frak A}\times{\frak A}\times{\frak A}\times{\frak A}^\ast\to {\mathbb C}$ we construct the ternary algebra $({\frak A},\tau)$ by defining the ternary multiplication $\tau$ as follows:
\begin{equation}
\theta(\tau(a,b,c))=T(a,b,c,\theta)
\label{dual_T}
\end{equation}
where $a,b,c\in {\frak A},\;\theta\in{\frak A}^\ast$. Particularly given a $\mathbb C$-bilinear mapping $L:{\frak A}\times{\frak A}\to \mbox{Lin}({\frak A})$, where $\mbox{Lin}({\frak A})$ is the algebra of linear operators of a vector space $\frak A$, we define
\begin{equation}
T(a,b,c,\theta)=\theta(L(a,b)\cdot c)
\end{equation}
and applying (\ref{dual_T}) we get the ternary algebra $({\frak A},\tau)$ whose ternary multiplication $\tau$ can be described implicitly by the formula
\begin{equation}
\tau(a,b,c)=L(a,b)\cdot c, \quad a,b,c\in {\frak A}
\label{operator_way}
\end{equation}
Applying this construction to a module over an associative unital algebra we can construct a ternary algebra by means of (\ref{dual_T}) or (\ref{operator_way}). Indeed if $\frak A$ is a left $\cal A$-module, where $\cal A$ is a binary unital associative complex algebra, ${\frak A}^\ast$ is the dual module and $T:{\frak A}\times{\frak A}\times{\frak A}\times{\frak A}^\ast\to {\cal A}$ is an $\cal A$-multilinear mapping  then $({\frak A},\tau)$ is the ternary algebra with the ternary multiplication defined by (\ref{dual_T}). Similarly given an $\cal A$-module $\frak A$ and a $\cal A$-bilinear mapping $L:{\frak A}\times{\frak A}\to \mbox{Lin}({\frak A})$, where $\mbox{Lin}({\frak A})$ is the algebra of $\cal A$-linear operators of a module $\frak A$, then ${\frak A}$ is the ternary algebra $({\frak A},\tau)$ with the ternary multiplication $\tau$ defined by (\ref{operator_way}).

We can use (\ref{dual_T},\ref{operator_way}) to construct the ternary algebras by means of well known structures of differential geometry on a smooth manifold. Let $\pi:E\to M$ be a vector bundle over a smooth finite dimensional manifold $M$, $C^{\infty}(M)$ be the algebra of smooth functions on a smooth manifold $M$, and ${\Gamma}(E)$ be the module of smooth sections of $E$. Given a $C^{\infty}(M)$-multilinear mapping $T:{\Gamma}(E)\times{\Gamma}(E)\times{\Gamma}(E)\times{\Gamma}(E^{\ast})\to C^{\infty}(M)$, where $E^{\ast}$ is the dual bundle, we obtain the ternary algebra $(\Gamma(E),\tau)$ of smooth sections of a vector bundle $E$ with the ternary multiplication $\theta(\tau(\xi,\eta,\chi))=T(\xi,\eta,\chi,\theta)$, where $\xi,\eta,\chi$ are sections of $E$, and $\theta\in \Gamma(E^{\ast})$. Particularly if $M$ is a smooth manifold, $E=TM$ is the tangent bundle, $E^{\ast}=T^{\ast}M$ is the cotangent bundle, $\Gamma(E)={\frak D}(M)$ is the module of vector fields, $\Gamma(E^{\ast})=\Omega^1(M)$ is the module of 1-forms, $\nabla$ is an affine connection on $M$, and
\begin{equation}
R(X,Y)=\nabla_X\,\nabla_Y-\nabla_Y\,\nabla_X-\nabla_{[X,Y]}, \quad X,Y\in {\frak D}(M)
\end{equation}
is the curvature of $\nabla$ then we have the $C^{\infty}(M)$-multilinear mapping
\begin{equation}
T:{\frak D}(M)\times{\frak D}(M)\times{\frak D}(M)\times\Omega^1(M)\to C^{\infty}(M)
\end{equation}
induced by the curvature
\begin{equation}
T(X,Y,Z,\omega)=\omega\,(R(X,Y)\cdot Z)
\end{equation}
and this mapping induces the structure of the ternary algebra $({\frak D}(M),\tau)$ on the module of vector fields with the ternary multiplication
\begin{equation}
\tau(X,Y,Z)=R(X,Y)\cdot Z
\end{equation}

Let ${\cal A}, {\cal B}$ be (binary) unital associative algebras over $\mathbb C$ with involutions respectively $a\to a^{\star}$ and $b\to b^{\ast}$, where $a\in {\cal A}, b\in {\cal B}$. Let $\cal M$ be a ${\cal A}-{\cal B}$-bimodule. We suppose that $\overline{\cal M}$ is an Abelian group which is isomorphic to the Abelian group $\cal M$, where $\varepsilon: m\in {\cal M}\to \varepsilon(m)={\bar m}\in \overline{\cal M}$ is the corresponding isomorphism. Then $\overline{\cal M}$ can be endowed with the structure of ${\cal B}-{\cal A}$-bimodule if we define the right and left multiplication by elements of algebras ${\cal A}$, ${\cal B}$ as follows
\begin{equation}
\varepsilon(m)\cdot a=\varepsilon(a^{\star}\cdot m),\quad
  b\cdot \varepsilon(m)=\varepsilon({m\cdot b^{\ast}}),\quad
       \forall\, m\in{\cal M},\; a\in {\cal A},\; b\in {\cal B}
\end{equation}
Let ${\cal M}\otimes_{\cal B}\overline{\cal M}$, $\overline{\cal M}\otimes_{\cal A}{\cal M}$ be the tensor products of modules where the first tensor product has the structure of $\cal A$-bimodule, and the second has the structure of $\cal B$-bimodule. It is clear that the algebras ${\cal A}$, $\cal B$ can be viewed respectively as $\cal A$-bimodule and $\cal B$-bimodule. We also assume that there are two homomorphisms
$\varphi\!: {\cal M}\otimes_{\cal B}\overline{\cal M}\to {\cal A}$, $\psi\!: \overline{\cal M}\otimes_{\cal A}{\cal M}\to {\cal B}$ respectively of $\cal A$-bimodules and $\cal B$-bimodules which satisfy
\begin{equation}
(\varphi(m\otimes {\bar n}))^{\star}=\varphi(n\otimes {\bar m}),\;\;
   (\psi({\bar m}\otimes {n}))^{\ast}=\psi({\bar n}\otimes {m}),\;\;
         \varphi(m\otimes {\bar n})\cdot p=m\cdot \psi ({\bar n}\otimes p)
\end{equation}
where $m,n,p\in {\cal M}$.
Evidently $\cal M$ has the structure of vector space over $\mathbb C$. We define the ternary law of composition $\tau$ on ${\cal M}$ by the formula
\begin{equation}
\tau(m,n,p)=\varphi(m\otimes {\bar n})\cdot p,\quad m,n,p\in {\cal M}
\end{equation}
\begin{proposition}
$({\cal M}, \tau)$ is {\alr}-partially associative ternary algebra.
\label{lr-associativity-of-M}
\end{proposition}
\noindent
Indeed for any quintuple of elements $m,n,p,q,r$ of $\cal M$ we have
\begin{equation}
\tau(\tau(m,n,p),q.r)=\tau(\varphi(m\otimes {\bar n})\cdot p,q,r)
       =\varphi\big((\varphi(m\otimes {\bar n})\cdot p)\otimes {\bar q}\big)\cdot r
       =\big(\varphi(m\otimes {\bar n})\varphi(p\otimes {\bar q})\big)\cdot r
\end{equation}
On the other hand
\begin{equation}
\tau(m,n,\tau(p,q,r))=\big(\varphi(m\otimes {\bar n})\varphi(p\otimes {\bar q})\big)\cdot r
\end{equation}
and this ends the proof.

From this proposition and the Proposition \ref{lr-associativity}
it follows that for any $n\in {\cal M}$ the binary algebra $({\cal M},\tau_{n,2})$, where $\tau_{n,2}(m,p)=\varphi(m\otimes{\bar n})\cdot p$, is an associative algebra.

Let us define the ternary multiplication $\sigma$ on ${\cal M}$ by $\sigma(m,n,p)=\varphi(m\otimes {\bar n})\cdot p+m\cdot \psi({\bar n}\otimes p)$, where $m,n,p\in {\cal M}$.
\begin{proposition}
The ternary algebra $({\cal M},\sigma)$ is the ternary algebra of Jordan type, and the ternary multiplication $\sigma$ of this algebra satisfies the identity
\begin{equation}
\sigma(m,n,\sigma(p,q,r))- \sigma(p,q,\sigma(m,n,r))+\sigma(\sigma(p,q,m),n,r)-\sigma(m,\sigma(q,p,n),r)=0
\label{identity}
\end{equation}
where $m,n,p,q,r\in {\cal M}$.
\end{proposition}
We use the ternary algebra $({\cal M}, \sigma)$ to construct a Lie algebra which will be constructed by means of defining commutation relations. We denote the set of generators of this Lie algebra by ${U}_{\bar n}, S_{pq}, U_m$, where $m,n,p,q\in {\cal M}$, i.e. we assign to each element $m\in {\cal M}$ the generator $U_m$, to each element ${\bar n}\in \bar{\cal M}$ the generator $U_{\bar n}$, and to each pair $(p,q)\in {\cal M}\times {\cal M}$ the generator $S_{pq}$. Let $\{U_m\}, \{U_{\bar n}\}, \{S_{pq}\}$ be the linear spans induced by the corresponding generators and ${\cal L}=\{U_m\}\oplus\{U_{\bar n}\}\oplus \{S_{pq}\}$. We define
\begin{gather}
[U_m,U_{\bar n}] = S_{mn},\quad\;\;\;\;\;\;\;\;
  [S_{pq},U_{m}] = U_{\sigma(p,q,m)}
\label{brackets-1}\\
    [S_{pq},U_{\bar n}] = -U_{\overline{\sigma(q,p,n)}},\;\;
   [S_{mn},S_{pq}] = S_{\sigma(m,n,p)\,q}-S_{p\,\sigma(n,m,q)}
\label{brackets-2}
\end{gather}
\begin{proposition}
The vector space $\cal L$ endowed with the brackets defined by (\ref{brackets-1}), (\ref{brackets-2}) is the Lie algebra, and the identity (\ref{identity}) leads to the Jacoby identity for the brackets (\ref{brackets-1}), (\ref{brackets-2}).
\end{proposition}
We can construct a matrix representation for the Lie algebra $\cal L$ if we consider the set $\mbox{Mat}_2({\cal A},{\cal B},{\cal M})$ of all $2\times 2$ square matrices of the type
\begin{equation}
A=\left(
  \begin{array}{cc}
    a & m \\
    {\bar n} & b \\
  \end{array}
\right), \quad
a\in {\cal A},\; b\in {\cal B},\; m\in {\cal M},\; {\bar n}\in \bar{\cal M}
\end{equation}
Given two such matrices
\begin{equation}
A=\left(
  \begin{array}{cc}
    a & m \\
    {\bar n} & b \\
  \end{array}
\right), \quad
B=\left(
  \begin{array}{cc}
    a^\prime & m^\prime \\
    {\bar n}^\prime & b^\prime \\
  \end{array}
\right)
\end{equation}
we define their product as follows
\begin{equation}
AB=\left(
  \begin{array}{cc}
    a\,a^\prime+\varphi(m\otimes {\bar n}^\prime)& a\cdot m^\prime+m\cdot b^\prime\\
    {\bar n}\cdot a^\prime+b\cdot {\bar n}^\prime & b\,b^\prime+\psi({\bar n}\otimes m^\prime) \\
  \end{array}
\right)
\label{product of matrices}
\end{equation}
\begin{proposition}
The vector space of matrices $\mbox{Mat}_2({\cal A},{\cal B},{\cal M})$ endowed with the multiplication (\ref{product of matrices}) is a unital associative (binary) algebra with the unity element
\begin{eqnarray}
E=\left(
  \begin{array}{cc}
    e & 0 \\
    {\bar 0} & e^\prime \\
  \end{array}
\right)
\end{eqnarray}
where $e$ is the unity element of $\cal A$ and $e^\prime$ is the unity element of $\cal B$.
\end{proposition}
Now we construct the matrix representation for the Lie algebra $\cal L$ as follows
\begin{equation}
U_{m}=\left(
  \begin{array}{cc}
    0 & m \\
    {\bar 0} & 0,\\
  \end{array}
\right),\quad
U_{\bar n}=\left(
  \begin{array}{cc}
    0 & 0 \\
    {\bar n} & 0 \\
  \end{array}
\right)
\end{equation}
and the matrix representation of $S_{pq}$ can be found by explicit calculation.

Having constructed the Lie algebra $\cal L$ and its matrix representation we can go further and construct a gauge field theory based on ternary algebra $({\cal M},\sigma)$. For this purpose we take a vector bundle $E$ over a smooth finite-dimensional manifold $M$ with the fiber $({\cal M},\sigma)$. A section of this bundle is a gauge field of our theory may be called ternon. Now in each fiber we construct the Lie algebra $\cal L$ and this leads us to the vector bundle of Lie algebras. Next we construct the Lie group corresponding to $\cal L$ by means of exponential mapping and Campbell-Hausdorf series. We get the principal fiber bundle and the we proceed in constructing the gauge field theory as usual.


\section{Associative multiplications of cubic matrices}

In this section we consider a vector space of cubic matrices, where by cubic matrix we mean a quantity $A=(A_{ijk})$ with three subscripts $i,j,k$ each running some set of integers. We use this vector space to construct a ternary algebra by means of triple product of cubic matrices. A triple product or ternary multiplication of cubic matrices is constructed in analogy with the classical product of two rectangular matrices by means of summation which is taken over certain system of subscripts of three cubic matrices. Our aim in this section is to find all totally associative ternary multiplications of first or second kind, and we prove that there are four ternary multiplications of cubic matrices each yielding the associative ternary algebra of second kind.

Let $A=(A_{k\hat{m}\bar{n}})$, where $A_{k\hat{m}\bar{n}}\in {\mathbb C}$ and $k,\hat{m},\bar{n}$ are integers satisfying $1\leq k\leq K$, $1\leq \hat{m}\leq M$, $1\leq \bar{n}\leq N$. We will call $A$ a complex $KMN$--\,space matrix provided that its entries $A_{k\hat{m}\bar{n}}$ are arranged in the vertices of  a 3--\,dimensional lattice and this structure is shown in particular case of a cubic matrix on the figure below. Let us denote the set of all such matrices by $\mbox{SMat}_{KMN}({\mathbb C})$, i.e.
\[
\mbox{SMat}_{KMN}({\mathbb C}) \!=\! \left\{ A \!= \!(A_{k\hat{m}\bar{n}}) \!:A_{k\hat{m}\bar{n}}\in {\mathbb C}, k=1,2,\ldots,K; \hat{m}=1,2,\ldots,M; \bar{n}=1,2,\ldots,N \right\}
\]
The set of $KMN$--\,space matrices is the vector space if we define the addition of space matrices and multiplication by complex numbers as usual
\begin{eqnarray}
A+B=(A_{k\hat{m}\bar{n}}+B_{k\hat{m}\bar{n}}),\quad
         \lambda\,A=(\lambda\,A_{k\hat{m}\bar{n}}),\quad \lambda\in {\mathbb C}
\end{eqnarray}
Our main concern in this paper is a special case of space matrices when $K=M=N$. In this case we will call $A=(A_{ijk})$, where $i,j,k=1,2,\ldots,N$, a complex $N$-cubic matrix and denote the vector space of such matrices by $\mbox{CMat}_{N}({\mathbb C})$.
Particularly, if $A=(A_{ijk})\in \mbox{CMat}_{3}({\mathbb C})$ is a cubic matrix of third order then we will place its entries into the vertices of $3$--dimensional lattice as follows
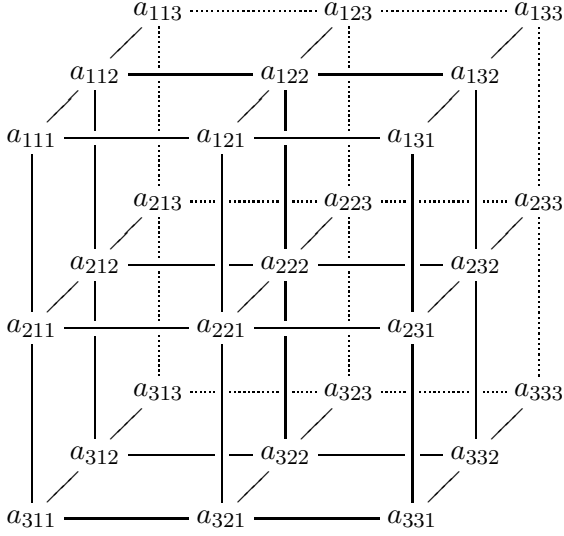
\begin{eqnarray}
\xymatrix@!0{
& & a_{113}   \ar@{.}[rrr]\ar@{.}'[d]'[dd][ddd]
& & & a_{123} \ar@{.}[rrr]\ar@{.}'[d]'[dd][ddd]
& & & a_{133} \ar@{.}[ddd]
\\
& a_{112}     \ar@{-}[ur] \ar@{-}[rrr] \ar@{-}'[d][ddd]
& & & a_{122} \ar@{-}[ur] \ar@{-}[rrr] \ar@{-}'[d][ddd]
& & & a_{132} \ar@{-}[ur]              \ar@{-}[ddd]
\\
a_{111}       \ar@{-}[ur] \ar@{-}[rrr] \ar@{-}[ddd]
& & & a_{121} \ar@{-}[ur] \ar@{-}[rrr] \ar@{-}[ddd]
& & & a_{131} \ar@{-}[ur] \ar@{-}[ddd] \ar@{-}[ddd]
\\
& & a_{213}   \ar@{.}'[r]'[rr][rrr] \ar@{.}'[d]'[dd][ddd]
& & & a_{223} \ar@{.}'[r]'[rr][rrr] \ar@{.}'[d]'[dd][ddd]
& & & a_{233} \ar@{.}[ddd]
\\
& a_{212}     \ar@{-}[ur] \ar@{-}'[rr][rrr] \ar@{-}'[d][ddd]
& & & a_{222} \ar@{-}[ur] \ar@{-}'[rr][rrr] \ar@{-}'[d][ddd]
& & & a_{232} \ar@{-}[ur] \ar@{-}[ddd]
\\
a_{211}       \ar@{-}[ur] \ar@{-}[rrr] \ar@{-}[ddd]
& & & a_{221} \ar@{-}[ur] \ar@{-}[rrr] \ar@{-}[ddd]
& & & a_{231} \ar@{-}[ur]              \ar@{-}[ddd]
\\
& & a_{313}   \ar@{.}'[r]'[rr][rrr]
& & & a_{323} \ar@{.}'[r]'[rr][rrr]
& & & a_{333}
\\
& a_{312}     \ar@{-}[ur] \ar@{-}'[rr][rrr]
& & & a_{322} \ar@{-}[ur] \ar@{-}'[rr][rrr]
& & & a_{332} \ar@{-}[ur]
\\
a_{311}       \ar@{-}[rrr]\ar@{-}[ur]
& & & a_{321} \ar@{-}[rrr]\ar@{-}[ur]
& & & a_{331} \ar@{-}[ur]
}\nonumber
\end{eqnarray}
The above $3$--dimensional lattice clearly shows that if one fixes a value of subscript $k$ in $A_{ijk}$ with $i,j$ ranging from 1 to 3 then the corresponding entries of cubic matrix $A$ form the square matrix of order $3$. Hence we get three square matrices of order $3$ which all together give us a cubic matrix of third order $A$. Therefore any cubic matrix $A$ of third order can be represented as the set of square matrices of order $3$ as follows
\begin{eqnarray}
A=(A_{ijk})= \left[\begin{matrix}\left(\begin{matrix}a_{111}&a_{121}&a_{131}\\a_{211}&a_{221}&a_{231}\\a_{311}&a_{321}&a_{331}\end{matrix}\right)_{k=1}&\left(\begin{matrix}a_{112}&a_{122}&a_{132}\\a_{212}&a_{222}&a_{232}\\a_{312}&a_{322}&a_{332}\end{matrix}\right)_{k=2}&\left(\begin{matrix}a_{113}&a_{123}&a_{133}\\a_{213}&a_{223}&a_{233}\\a_{313}&a_{323}&a_{333}\end{matrix}\right)_{k=3}\end{matrix}\right]
\nonumber
\end{eqnarray}


Now our aim is to construct a multiplication of space matrices. We will do this for cubic matrices of order $N$ because constructed multiplication can be extended to a vector space of $KMN$--\,space matrices in an obvious way. If $A=(A_{ijk})\in \mbox{CMat}_{N}({\mathbb C})$ is a cubic matrix of order $N$ then it induces two mappings as follows:
\begin{enumerate}
\itemsep-1pt
\item[i)]
$
A_{\mbox{vec}\to \mbox{op}}: x=(x_k)\in {\mathbb C}^N\to L=(L_{ij})=(\sum_k A_{ijk}\,x_k)\in \mbox{Lin}({\mathbb C}^N)
$
\item[ii)]
$
A_{\mbox{op}\to \mbox{vec}}: L=(L_{jk})\in \mbox{Lin}({\mathbb C}^N)\to x=(x_i)=(\sum_{j,k} A_{ijk}\,L_{jk})\in {\mathbb C}^N
$
\end{enumerate}
where $\mbox{Lin}({\mathbb C}^N)$ is the vector space of linear operators acting on ${\mathbb C}^N$.
These mappings determined by a cubic matrix of order $N$ show that if we wish to construct a multiplication of cubic matrices of order $N$ which is based on composition of linear mappings as in the calculus of rectangular matrices  then we should take three cubic matrices in order to close a corresponding multiplication in the sense that the product of three cubic matrices of order $N$ will be the cubic matrix of order $N$. Given two cubic matrices $A,B\in \mbox{CMat}_{N}({\mathbb C})$ we can form two products of corresponding mappings $(A_{\mbox{vec}\to \mbox{op}})\circ (B_{\mbox{op}\to \mbox{vec}})$ and $(A_{\mbox{op}\to \mbox{vec}})\circ (B_{\mbox{vec}\to \mbox{op}})$ which are not close with respect to composition because the first product is the linear mapping ${\mathbb C}^N\to {\mathbb C}^N$, and the second is the linear mapping $\mbox{Lin}({\mathbb C}^N)\to \mbox{Lin}({\mathbb C}^N)$. It is obvious than we can close a procedure of taking compositions of this sort of mappings by adding one more cubic matrix $C=(C_{ijk})\in \mbox{CMat}_{N}({\mathbb C})$. In this case the triple product $(A_{\mbox{vec}\to \mbox{op}})\circ (B_{\mbox{op}\to \mbox{vec}})\circ (C_{{\mbox{vec}\to \mbox{op}}})$ closes the operation of taking compositions giving the mapping ${\mathbb C}^N\to \mbox{Lin}({\mathbb C}^N)$.
This kind of reasoning suggests a possible way of constructing the cubic matrix of order $N$ from three given ones by means of summation which is taken over the certain pair of subscripts.

Given three cubic matrices of order $N$ we have nine subscripts, and because the number of possible combinations  of subscripts is finite we can use the methods of computer algebra to find all associative ternary multiplications either of first or second kind. Our analysis shows that there is no total associative ternary multiplication of first kind, and all total associative ternary multiplications of second kind are described by the following theorem.

\begin{theorem}
There are only four different triple products of complex cubic matrices of order $N$ which obey the total ternary associativity of second kind. These are
\begin{enumerate}
\itemsep-3pt
\item[1)]
${(A \odot B \odot C)}_{ijk} =\sum_{l,m,n} A_{ilm}B_{nlm}C_{njk},\quad
     {A \odot B \odot C} \rightarrow
\xy <1cm,0cm>:
(1,0)*+{A} , (2,0)*+{B} , (3,0)*+{C} ,
(1.25,-0.2)*+{\bullet} , (1.45,-0.2)*+{\circ} , (1.65,-0.2)*+{\circ} ,
(2.25,-0.2)*+{\circ} , (2.45,-0.2)*+{\circ} , (2.65,-0.2)*+{\circ} ,
(3.25,-0.2)*+{\circ} , (3.45,-0.2)*+{\bullet} , (3.65,-0.2)*+{\bullet} ,
(1.45,-0.26);(2.45,-0.26)**\crv{(1.5, -0.5)&(1.95, -0.7)&(2.4, -0.5)} ,
(1.65,-0.26);(2.65,-0.26)**\crv{(1.7, -0.4)&(2.15, -0.6)&(2.6, -0.4)} ,
(2.25,-0.26);(3.25,-0.26)**\crv{(2.3, -0.5)&(2.75, -0.7)&(3.2, -0.5)}
\endxy$
\item[2)]
${(A \odot B \odot C)}_{ijk} =\sum_{l,m,n} A_{ilm}B_{nml}C_{njk},\quad
     {A \odot B \odot C} \rightarrow
\xy <1cm,0cm>:
(1,0)*+{A} , (2,0)*+{B} , (3,0)*+{C} ,
(1.25,-0.2)*+{\bullet} , (1.45,-0.2)*+{\circ} , (1.65,-0.2)*+{\circ} ,
(2.25,-0.2)*+{\circ} , (2.45,-0.2)*+{\circ} , (2.65,-0.2)*+{\circ} ,
(3.25,-0.2)*+{\circ} , (3.45,-0.2)*+{\bullet} , (3.65,-0.2)*+{\bullet} ,
(1.45,-0.26);(2.65,-0.26)**\crv{(1.5, -0.5)&(2.05, -0.7)&(2.6, -0.5)} ,
(1.65,-0.26);(2.45,-0.26)**\crv{(1.7, -0.4)&(2.05, -0.6)&(2.4, -0.4)} ,
(2.25,-0.26);(3.25,-0.26)**\crv{(2.3, -0.5)&(2.75, -0.7)&(3.2, -0.5)}
\endxy$
\item[3)]
${(A \odot B \odot C)}_{ijk} =\sum_{l,m,n} A_{ijl}B_{nml}C_{mnk},\quad
    {A \odot B \odot C} \rightarrow
\xy <1cm,0cm>:
(1,0)*+{A} , (2,0)*+{B} , (3,0)*+{C} ,
(1.25,-0.2)*+{\bullet} , (1.45,-0.2)*+{\bullet} , (1.65,-0.2)*+{\circ} ,
(2.25,-0.2)*+{\circ} , (2.45,-0.2)*+{\circ} , (2.65,-0.2)*+{\circ} ,
(3.25,-0.2)*+{\circ} , (3.45,-0.2)*+{\circ} , (3.65,-0.2)*+{\bullet} ,
(1.65,-0.26);(2.65,-0.26)**\crv{(1.7, -0.5)&(2.15, -0.7)&(2.6, -0.5)} ,
(2.25,-0.26);(3.45,-0.26)**\crv{(2.3, -0.5)&(2.85, -0.7)&(3.4, -0.5)} ,
(2.45,-0.26);(3.25,-0.26)**\crv{(2.45, -0.4)&(2.85, -0.6)&(3.2, -0.4)} ,
\endxy$
\item[4)]
${(A \odot B \odot C)}_{ijk} =\sum_{l,m,n} A_{ijl}B_{mnl}C_{mnk},\quad
      {A \odot B \odot C} \rightarrow
\xy <1cm,0cm>:
(1,0)*+{A} , (2,0)*+{B} , (3,0)*+{C} ,
(1.25,-0.2)*+{\bullet} , (1.45,-0.2)*+{\bullet} , (1.65,-0.2)*+{\circ} ,
(2.25,-0.2)*+{\circ} , (2.45,-0.2)*+{\circ} , (2.65,-0.2)*+{\circ} ,
(3.25,-0.2)*+{\circ} , (3.45,-0.2)*+{\circ} , (3.65,-0.2)*+{\bullet} ,
(1.65,-0.26);(2.65,-0.26)**\crv{(1.7, -0.5)&(2.15, -0.7)&(2.6, -0.5)} ,
(2.25,-0.26);(3.25,-0.26)**\crv{(2.3, -0.5)&(2.75, -0.7)&(3.2, -0.5)} ,
(2.45,-0.26);(3.45,-0.26)**\crv{(2.45, -0.4)&(2.95, -0.6)&(3.4, -0.4)} ,
\endxy$
\label{theorem}
\end{enumerate}
\end{theorem}


\section{Ternary analog of Lie algebra}

Now our aim is to construct a ternary analog of Lie algebra with the help of ternary multiplication of cubic matrices and ternary analog of Lie bracket. A notion of ternary analog of Lie bracket is based on a faithful representation  of cyclic group $\mathbb{Z}_3$ by cubic roots of unity. We construct a ternary analog of Lie algebra which may be considered as an analog of the Lie algebra generated by Pauli matrices. In our construction we use the cubic matrices of second order with certain symmetries with respect to the subscripts as generators of our algebra and find all commutation relations. In this section we use the ternary multiplication of cubic matrices which has the property that any cyclic permutation of the matrices in the product is equivalent to the same permutation on the subscripts and this multiplication is studied in \cite{Abr-Kerner-Roy,Kerner_1,Kerner_2,Kerner_3, Vainerman-Kerner}.

Let $A,$ $B,$ $C$ be cubic matrices of order $N$. We define the triple product $A \circledcirc B \circledcirc C$ by the following formula
\begin{equation}
{(A \circledcirc B \circledcirc C)}_{ikl}
=\sum_{p,q,r} A_{piq}B_{qkr}C_{rlp},\quad
{A \circledcirc B \circledcirc C} \rightarrow
\xy <1cm,0cm>:
(1,0)*+{A} , (2,0)*+{B} , (3,0)*+{C} ,
(1.25,-0.2)*+{\circ} , (1.45,-0.2)*+{\bullet} , (1.65,-0.2)*+{\circ} ,
(2.25,-0.2)*+{\circ} , (2.45,-0.2)*+{\bullet} , (2.65,-0.2)*+{\circ} ,
(3.25,-0.2)*+{\circ} , (3.45,-0.2)*+{\bullet} , (3.65,-0.2)*+{\circ} ,
(1.25,-0.26);(3.65,-0.26)**\crv{(1.40, -0.5)&(2.45, -0.65)&(3.50, -0.5)} ,
(1.65,-0.26);(2.25,-0.26)**\crv{(1.75, -0.38)&(1.95, -0.42)&(2.15, -0.38)},
(2.65,-0.26);(3.25,-0.26)**\crv{(2.75, -0.38)&(2.95, -0.42)&(3.15, -0.38)}
\endxy
\label{multiplication}
\end{equation}
It is easy to see that any cyclic permutation of the matrices in the product is equivalent to the same permutation on the subscripts, i.e.
${(A \circledcirc B \circledcirc C)}_{ikl} = {(B \circledcirc C \circledcirc A)}_{kli} = {(C \circledcirc A \circledcirc B)}_{lik}.$  It should be mentioned that the ternary multiplication (\ref{multiplication}) is neither partially nor totally associative.

The ternary algebra of cubic matrices of order $N$ with respect to the multiplication (\ref{multiplication}) can be decomposed into direct sum of subspaces of cubic matrices with certain symmetries according to the irreducible representation of the symmetry group $S_3.$ It should be mentioned that we have the similar decomposition in the case of square matrices of order $N,$ where the algebra of square matrices can be decomposed into the direct sum of subspaces of symmetric and skew-symmetric matrices according to representation of the permutation group $S_2.$
The symmetry group $S_3$ possesses a full and faithful representation on the complex plane, which can be generated by two elements representing a cyclic and an odd permutation. This representation can be constructed by assigning the operator of multiplication by the cubic root of unity $j=e^{\frac{2\pi i}{3}}$ to the cyclic permutation $\left ( \begin{smallmatrix} abc \\ bca \end{smallmatrix} \right )$ and assigning the operator of complex conjugation to the odd permutation  $\left ( \begin{smallmatrix} abc \\ cab \end{smallmatrix} \right )$. In order to construct the subspace of cubic matrices with certain symmetries we use the cyclic part of this representation where the cyclic group $\mathbb{Z}_3$ is represented by cubic roots of unity.

A cubic matrix $\rho=(\rho_{ikl})$ of order $N$ is said to be a $j$-skew-symmetric if it satisfies
\begin{eqnarray}
\rho_{ikl}=j\rho_{kli}=j^2\rho_{lik}
\nonumber
\end{eqnarray}
Similarly  a cubic matrix $\overline{\rho}=(\overline{\rho}_{ikl})$ of order $N$ is said to be a $j^2$-skew-symmetric if it satisfies
\begin{eqnarray}
\overline{\rho}_{ikl}=j\overline{\rho}_{kli}=j^2\overline{\rho}_{lik}
\nonumber
\end{eqnarray}
It can be shown that the space of cubic matrices of order $N$ can be decomposed into the direct sum of the subspace of $j$-skew-symmetric matrices, the subspace of $j^2$-skew-symmetric matrices and the subspace of symmetric matrices, where under symmetric matrix we mean a cubic matrix $\omega=(\omega_{ikl})$ of order $N$ which satisfies $\omega_{ijk}=\omega_{jki}=\omega_{kji}$.

Each of subspaces of $j$-skew-symmetric matrices and of $j^2$-skew-symmetric matrices has the dimension ${(N^3-N)}/{3}.$ The subspace of symmetric matrices can be decomposed into the direct sum of the subspace of diagonal cubic matrices (all entries excepting diagonal ones are equal to zero) and the subspace of symmetric matrices with zeros on diagonal. Obviously the dimension of the subspace of diagonal cubic matrices is equal to $N$ and the subspace of symmetric matrices with zeros on diagonal is equal to ${(N^3-N)}/{3}.$ Let $\{\rho^{\a}\}$ be a basis for the subspace of $j$-skew-symmetric matrices, $\{\overline{\rho}^{\a}\}$ be a basis for the subspace of $j^2$-skew-symmetric matrices, and $\{\omega^{\a}\}$ be a basis for the subspace of symmetric matrices with zeros on diagonal, where $\a=1,2,\ldots ,{(N^3-N)}/{3}.$ A basis for the subspace of diagonal cubic matrices will be denoted by $\{\eta^{\a}\},$ where $\a=1,2,\ldots,N.$ Thus we have
\begin{eqnarray}
\lefteqn{\rho_{ikl}^{\a}=j\rho_{kli}^{\a}=j^2\rho_{lik}^{\a},\quad \a =1,2,\ldots, \frac{N^3-N}{3}}\nonumber\\
\lefteqn{\overline{\rho}^{\,\a}_{ikl}=j^2\overline{\rho}^{\,\a}_{kli}=j\overline{\rho}^{\,\a}_{lik} \quad \a=1,2,\ldots ,\frac{N^3-N}{3}},\nonumber\\
\lefteqn{\omega_{ijk}^{\a}=\omega_{jki}^{\a}=\omega_{kji}^{\a}, \quad \omega^{\a}_{iii}=0, \quad \a=1,2,\ldots ,\frac{N^3-N}{3}}.\nonumber
\end{eqnarray}

 Now our aim is to construct a ternary analog of Lie algebra by means of the ternary $j$-commutator which we define as follows
 \begin{equation}
 [A,B,C]=A\circledcirc B\circledcirc C+jB\circledcirc C\circledcirc A+j^2C\circledcirc A\circledcirc B,\quad A,B,C\in \mbox{CMat}_{N}({\mathbb C})
 \label{commutator}
 \end{equation}
 where the ternary multiplication $A\circledcirc B\circledcirc C$ is defined by (\ref{multiplication}). It is easy to see that the ternary $j$-commutator (\ref{commutator}) may be viewed as an analog of the binary Lie commutator $[X,Y]=XY-YX,$ where $X,$ $Y$ are the elements of an associative binary algebra, because we replace the symmetric group $\mathbb{Z}_2$ and its representation $\{1,-1\}$ respectively by the cyclic group $\mathbb{Z}_3$ and its representation $\{1,j,j^2\}.$ Indeed in the case of the binary Lie commutator we have the property $[X,Y]=-[Y,X].$ For the ternary $j$-commutator we have the similar property $[A,B,C]=j[B,C,A]=j^2[C,A,B].$  In analogy with the binary Lie commutator for any $A\in \mbox{CMat}_{N}({\mathbb C})$ we have $[A,A,A]=0.$

Let us consider the space of cubic matrices of order $2,$ i.e. each subscript runs from 1 to 2.  In this case all subspaces mentioned above are 2-dimensional, which means that basis for each subspace consists of two matrices, i.e. $\a=1,2.$  We construct the basis for the subspace of diagonal cubic matrices $\{\eta^1,\eta^2 \}$ by choosing $\eta_{111}^1=1, \;\;\eta^1_{222}=0,\;\;\eta^2_{111}=0,\;\;\eta^2_{222}=1.$ We choose the basis for the subspace of symmetric cubic matrices with zeros on diagonal by fixing $\omega^1_{212}=\omega^1_{221}=\omega^1_{122}=1,$ $\omega^2_{121}=\omega^2_{112}=\omega^2_{211}=1$ with other entries are equal to zero. We get the basis for the subspace of $j$-skew-symmetric cubic matrices by taking $\rho_{212}^1=j\rho_{122}^1=j^2\rho_{221}^1=1,$ $\rho_{121}^2=j\rho_{211}^2=j^2\rho_{112}^2=1$ with other entries equal to zero and the basis for the subspace of $j^2$-skew-symmetric cubic matrices if we put $\orhoA_{212}=j^2\orhoA_{122}=j\orhoA_{221}=1,$ $\orhoB_{121}=j\orhoB_{211}=j^2\orhoB_{112}=1$ with other entries equal to zero. The space of cubic matrices of order 2 equipped with the ternary $j$-commutator can be considered as a ternary analog of matrix Lie algebra and all commutation relations of this algebra in the basis $\{\eta^1,\eta^2,\omega^1,\omega^2,\rho^1,\rho^2,\overline{\rho}^1,\overline{\rho}^2\}$ are given in the following table:

\begin{align*}
[\rhoA, \rhoB, \rhoA]       &= -\orhoB, & [\rhoB, \rhoA, \rhoB]        &= -\orhoA  &  [\orhoA, \orhoB, \orhoA] &= 2\orhoB, & [\orhoB, \orhoA, \orhoB] &= 2\orhoA \\
[\omegaA, \omegaB, \omegaA] &= -\orhoB, & [\omegaB, \omegaA, \omegaB]  &= -\orhoA  &  [\etaA,  \etaB,  \etaA]  &= 0,       & [\etaB,  \etaA,  \etaB]  &= 0       \\
\\
[\rhoA, \orhoA, \rhoA]      &= 3\orhoA, &  [\rhoB, \orhoB, \rhoB]      &= 3\orhoB, & [\orhoA, \rhoA, \orhoA]   &= 0,       & [\orhoB, \rhoB, \orhoB] &= 0        \\
[\rhoA, \orhoB, \rhoA]      &= -\orhoB, &  [\rhoB, \orhoA, \rhoB]      &= -\orhoA, & [\orhoA, \rhoB, \orhoA]   &= -\orhoB, & [\orhoB, \rhoA, \orhoB] &= -\orhoA  \\
[\rhoA, \orhoA, \rhoB]      &= 2\orhoB, &  [\rhoA, \orhoB, \rhoB]      &= 2\orhoA, & [\rhoA, \orhoA, \orhoB]   &= -\orhoB, & [\rhoB, \orhoB, \orhoA] &= -\orhoA  \\
\\
[\rhoA, \omegaA, \rhoA] &= 0,        &  [\omegaA, \rhoA, \omegaA] &= 3\orhoA, & [\rhoB, \omegaB, \rhoB]   &= 0,       & [\omegaB, \rhoB, \omegaB] &= 3\orhoB \\
[\rhoA, \omegaB, \rhoA] &= 2\orhoB,  &  [\omegaB, \rhoA, \omegaB] &= 2\orhoA, & [\rhoB, \omegaA, \rhoB]   &= 2\orhoA, & [\omegaA, \rhoB, \omegaA] &= 2\orhoB \\
[\rhoB, \rhoA, \omegaA] &= -\orhoB,  &  [\omegaB, \rhoB, \rhoA]   &= -\orhoA, & [\omegaA, \rhoA, \omegaB] &= -\orhoB, & [\omegaA, \rhoB, \omegaB] &= -\orhoA \\
\\
[\orhoA, \omegaA, \orhoA] &= 3\orhoA,  &  [\omegaA, \orhoA, \omegaA] &= 0,       & [\orhoB, \omegaB, \orhoB]  &= 3\orhoB, & [\omegaB, \orhoB, \omegaB] &= 0 \\
[\orhoA, \omegaB, \orhoA] &= -\orhoB,  &  [\omegaB, \orhoA, \omegaB] &= -\orhoA, & [\orhoB, \omegaA, \orhoB]  &= -\orhoA, & [\omegaA, \orhoB, \omegaA] &= -\orhoB \\
[\orhoA, \omegaA, \orhoB] &= -\orhoB,  &  [\orhoA, \omegaB, \orhoB]  &= -\orhoA, & [\orhoA, \omegaA, \omegaB] &= 2\orhoB, & [\omegaA, \omegaB, \orhoB] &= 2\orhoA \\
\\
[\rhoA, \etaA, \rhoA] &= \orhoA,  &  [\etaA, \rhoA, \etaA] &= 0,      & [\rhoB, \etaB, \rhoB]  &= \orhoB, & [\etaB, \rhoB, \etaB] &= 0 \\
[\rhoA, \etaB, \rhoA] &= \orhoB,  &  [\etaA, \rhoB, \etaA] &= \orhoB, & [\rhoB, \etaA, \rhoB]  &= \orhoA, & [\etaB, \rhoA, \etaB] &= \orhoA \\
[\rhoA, \rhoB, \etaA] &= \orhoB,  &  [\etaB, \rhoA, \etaB] &= \orhoA, & [\etaA, \rhoA, \etaB]  &= 0,      & [\etaA, \rhoB, \etaB] &= 0 \\
\\
[\orhoA, \etaA, \orhoA] &= \orhoA,  &  [\etaA,  \orhoA, \etaA]  &= 0,      & [\orhoB, \etaB,  \orhoB]  &= \orhoB, & [\etaB, \orhoB, \etaB] &= 0 \\
[\orhoA, \etaB, \orhoA] &= \orhoB,  &  [\etaA,  \orhoB, \etaA]  &= \orhoB, & [\orhoB, \etaA,  \orhoB]  &= \orhoA, & [\etaB, \orhoA, \etaB] &= \orhoA \\
[\orhoA, \etaA, \orhoB] &= \orhoB,  &  [\orhoA, \etaB,  \orhoB] &= \orhoA, & [\etaA,  \orhoA, \etaB]   &= 0,      & [\etaA, \orhoB, \etaB] &= 0 \\
\\
[\etaA,   \omegaA, \etaA]    &= 0,      &  [\etaA,   \omegaB, \etaA]    &= \orhoB, & [\etaB,   \omegaA, \etaB]   &= \orhoA,  & [\etaB,   \omegaB, \etaB]    &= \orhoB \\
[\omegaA, \etaA,   \omegaA]  &= \orhoA, &  [\omegaA, \etaB,   \omegaA]  &= \orhoB, & [\omegaB, \etaA,   \omegaB] &= \orhoA,  & [\omegaB, \etaB,   \omegaB]  &= \orhoB \\
[\etaA,   \omegaA, \etaB]    &= 0,       & [\etaA,   \omegaB, \etaB]    &= 0,       & [\etaA,   \omegaA, \omegaB] &= \orhoB,  & [\etaB,   \omegaB, \omegaA]  &= \orhoA \\
\\
[\rhoA,   \orhoA,  \omegaA]  &= 0,       & [\omegaA, \rhoA,   \orhoB]   &= 2\orhoB, & [\rhoA,   \orhoA,  \omegaB] &= -\orhoB, & [\rhoA,   \rhoB,   \omegaB]  &= -\orhoA \\
[\etaA,   \rhoA,   \orhoA]   &= \orhoA,  & [\etaA,   \rhoA,   \orhoB]   &= \orhoB,  & [\orhoA,  \etaB,   \rhoA]   &= \orhoB,  & [\rhoA,   \orhoB,  \etaB]    &= \orhoA \\
[\orhoA,  \omegaA, \rhoB]    &= -\orhoB, & [\rhoB,   \orhoB,  \omegaA]  &= -\orhoA, & [\omegaB, \rhoB,   \orhoA]  &= 2\orhoA, & [\rhoB,   \orhoB,  \omegaB]  &= 0 \\
[\etaA,   \orhoA,  \rhoB]    &= \orhoB,  & [\rhoB,   \etaA,   \orhoB]   &= \orhoA,  & [\etaB,   \rhoB,   \orhoA]  &= \orhoA,  & [\etaB,   \rhoB,   \orhoB]   &= \orhoB \\
[\etaA,   \orhoA,  \omegaA]  &= \orhoA,  & [\etaA,   \omegaB, \orhoA]   &= \orhoB,  & [\omegaA, \etaB,   \orhoA]  &= \orhoB,  & [\etaB,   \orhoA,  \omegaB]  &= \orhoA \\
[\etaA,   \orhoB,  \omegaA]  &= \orhoB,  & [\orhoB,  \etaA,   \omegaB]  &= \orhoA,  & [\orhoB,  \omegaA, \etaB]   &= \orhoA,  & [\etaB,   \orhoB,  \omegaB]  &= \orhoB \\
[\etaA,   \rhoA,   \omegaA]  &= \orhoA,  & [\rhoA,   \etaA,   \omegaB]  &= \orhoB,  & [\omegaA, \etaB,   \rhoA]   &= \orhoB,  & [\rhoA,   \etaB,   \omegaB]  &= \orhoA \\
[\omegaA, \etaA,   \rhoB]    &= \orhoB,  & [\omegaB, \etaA,   \rhoB]    &= \orhoA,  & [\omegaA, \etaB,   \rhoB]   &= \orhoA,  & [\etaB,   \rhoB,   \omegaB]  &= \orhoB
\end{align*}
\vskip.1cm
\noindent
We would like to point out that the above table of commutation relations demonstrates a peculiar property of the algebra of cubic matrices of second order in the basis $\{\eta^1,\eta^2,\omega^1,\omega^2,\rho^1,\rho^2,\overline{\rho}^1,\overline{\rho}^2\}$: the $j$-commutator of any three generators is proportional to one of the generators $\overline{\rho}^1,\overline{\rho}^2$. Moreover taking $j^2$-commutator
\begin{equation}
[A,B,C]^{*}=A\circledcirc B\circledcirc C+j^2\,B\circledcirc C\circledcirc A+j\,C\circledcirc A\circledcirc B,\quad A,B,C\in \mbox{CMat}_{N}({\mathbb C}),
\label{con-commutator}
\end{equation}
which can be considered as a conjugate commutator to $j$-commutator (\ref{commutator}), we get the similar table of commutation relations (we do not demonstrate it here in order not to overburden the paper) which clearly shows that the $j^2$-commutator of any three generators is proportional to one of the generators ${\rho}^1,{\rho}^2$. Finally considering a ternary analog of anti-commutator defined by
\begin{equation}
\{A,B,C\}=A\circledcirc B\circledcirc C+B\circledcirc C\circledcirc A+C\circledcirc A\circledcirc B,\quad A,B,C\in \mbox{CMat}_{N}({\mathbb C}),
\label{anti-commutator}
\end{equation}
and applying it to the cubic matrices of second order we obtain the table of commutation relations for the generators $\eta^1,$ $\eta^2,$ $\omega^1,$ $\omega^2,$ $\rho^1,$ $\rho^2,$ $\overline{\rho}^1$, $\overline{\rho}^2$ which can be shortly described as follows: the ternary anti-commutator of any three generators is a linear combination of the generators $\eta^1,\eta^2,\omega^1,\omega^2$.

We conclude this Section by pointing out an analogy between the generators of our  ternary algebra of cubic matrices of second order and the Pauli matrices
\begin{eqnarray}\sigma_1=\left ( \begin{array}{cc} 0 & 1 \\ 1 & 0 \end{array} \right ),\quad \sigma_2= \left ( \begin{array}{cc} 0 & -i \\ i & 0 \end{array} \right ),\quad \sigma_3=\left ( \begin{array}{cc} 1 & 0 \\ 0 & -1 \end{array} \right )
\nonumber
\end{eqnarray}
Indeed applying the ternary $j$-commutator to Pauli matrices we get
\begin{eqnarray} [i\sigma_1,i\sigma_2,i\sigma_1]=2i\sigma_2, \quad[i\sigma_2,i\sigma_1,i\sigma_2]=2i\sigma_1
\nonumber
\end{eqnarray}
In the table of commutation relations of cubic matrices of order 2 we have
\begin{eqnarray}
[ \orhoA, \orhoB, \orhoA ] = 2 \orhoB,\quad\;\;\;\; [\orhoB, \orhoA, \orhoB]=2\orhoA
\nonumber
\end{eqnarray}
which clearly demonstrates a striking analogy between the $j^2$-skew-symmetric generators of ternary algebra of cubic matrices of order 2 and the skew-Hermitian matrices $i\sigma_1,$ $i\sigma_2,$ $i\sigma_3,$ where $\sigma_1,$ $\sigma_2,$ $\sigma_3$ are Pauli matrices.

\section*{Acknowledgements}

V. Abramov, O. Liivapuu and S. Shitov gratefully acknowledge the financial support of their research by the Estonian Science Foundation under the research grant ETF 7427. V. Abramov, R. Kerner and O. Liivapuu would like to thank the French-Estonian Scientific Program "G. F. Parrot" for financial support under the travel grant P-2/2007.


\smallskip
\Date{Received January 16, 2008 
}

\label{lastpage}



\begin{thebibliography}{99}
\itemsep-3pt

{\small

\bibitem{Abramov-1}
V. Abramov. ${\mathbb Z}_3$-graded analogues of Clifford algebras and algebra of ${\mathbb Z}_3$-graded symmetries. Algebras, Groups and Geometries, {\bf 12} (1995), 201--221.

\bibitem{Abramov-2}
V. Abramov, Ternary generalizations of Grassmann algebra. Proc. Estonian Acad. Sci. Phys. Math. {\bf 45}  (1996), 152--160.

\bibitem{Abr-Kerner-Roy}
V. Abramov, R. Kerner, B. Le Roy. Hypersymmetry: a ${\mathbb Z}_3$-graded generalization of supersymmetry. J. Math. Phys. {\bf 38} (1997), 1650--1669.

\bibitem{Ataguema-thesis}
H. Ataguema. Classification et d\'{e}formations des alg\`{e}bres ternaires. PhD thesis. Universit\'{e} de Haute Alsace, 2008).

\bibitem{Ataguema-1}
H. Ataguema, A. Makhlouf. Notes on cohomologies of ternary algebras of associative type. Preprint arXiv: 0812.0707.

\bibitem{Ataguema-2}
H. Ataguema, A. Makhlouf, S. Silvestrov. Generalization of $n$-ary Nambu algebras and beyond. Preprint arXiv: 0812.4058.

\bibitem{Bars-Gunadyn-1}
I. Bars, M. G\"unaydin. Construction of Lie algebras and Lie superalgebras from ternary algebras. J. Math. Phys. {\bf 20} (1979), 1977--1993.

\bibitem{Bars-Gunadyn-2}
I. Bars, M. G\"unaydin. Dynamical theory of subconstituents based on ternary algebras. Phys. Rev. {\bf D22} (1980), 1403--1413.

\bibitem{Carlsson}
R. Carlsson. Cohomology of associative triple systems. Proc. American Math. Soc.
 {\bf 60} (1976), 1--7.

\bibitem{Kantor_1}
I. L. Kantor. Models of exceptional Lie algebras. Sov. Math. Dokl. {\bf 14} (1973), 254--258.

\bibitem{Kerner_1}
R. Kerner. Graduation $\mathbb Z_3$ et la racine cubique de l'\'{e}quation de Dirac. Comptes Rendus Acad. Sci. {\bf 312} (1991), 191--196.

\bibitem{Kerner_2}
R. Kerner. ${\mathbb Z}_3$-graded algebras and the cubic root of supersymmetry translations. J. Math. Phys. {\bf 33} (1992), 403--411.

\bibitem{Kerner_3}
R. Kerner. The cubic chessboard. Class. Quantum Grav. {\bf 14} (1997), A203-A225.

\bibitem{LeRoy_1}
B. Le Roy. A ${\mathbb Z}_3$-graded generalization of supermatrices. J. Math. Phys. {\bf 37} (1996), 474--483.

\bibitem{Vainerman-Kerner}
L. Vainerman, R. Kerner. On special classes of $n$-algebras. J. Math. Phys. {\bf 37} (1996), 2553-2565.

}

\end{thebibliography}
\end{document}